\documentclass{amsart}

\newtheorem{theorem}{Theorem}[section]
\newtheorem{lemma}[theorem]{Lemma}

\theoremstyle{definition}
\newtheorem{definition}[theorem]{Definition}

\theoremstyle{remark}

\newtheorem{proposition}[theorem]{Proposition}

\numberwithin{equation}{section}

\begin{document}

\title{AMS Journal Sample}

\author{Nazife Erkur\c{s}un \"Ozcan}
\address{ Department of Mathematics, Hacettepe University, Ankara, Turkey, 06800}
\email{erkursun.ozcan@hacettepe.edu.tr}

\subjclass{47A35, 47B42, 47B65, 47D03, 47D06}

\date{June 25, 2015}


\keywords{Markovian LR-sequences, KB-space, strong convergence, attractor}

\title{ Asymptotic behavior of operator sequences on KB-spaces}

\begin{abstract}
The concept of a constrictor was used by several mathematicians to characterize the asymptotic behavior of operators. In this paper we show that a positive LR-sequence on KB-spaces is mean ergodic if the LR-sequence has a weakly compact attractor. Moreover if the weakly compact attractor is an order interval, then a Markovian LR-sequence converges strongly to the finite dimensional fixed space. As a consequence we investigate also stability of LR-sequences of positive operators and existence of lower bound functions on KB-spaces. 
\end{abstract}

\maketitle

\section{Introduction}

 In \cite{ABEE}, it is proved that if $T$ is Markov operator on a KB-space then $T$ is mean ergodic and satisfies $\mathrm{dim (Fix} (T)) < + \infty$  whenever there exist a function $g \in E_+$ and a real $0 \leq \eta < 1$ such that $\displaystyle  \lim_{n\to\infty} \text{dist} ( \frac{1}{n} \sum_{k=0}^{n-1} T^k x, [-g,g]+\eta B_E ) =0$ for every $x \in E_+ \cap U_E$. In this paper, we extend this result to any Markov LR-sequence on KB-spaces.

Let $X$ be Banach space and $I_X = I$ be the identity operator. Let $E$ be a Banach lattice. Then $E_+ := \{ x\in E : x \geq 0 \}$ denotes the positive cone of $E$. On ${\mathcal L} (E)$ there is a canonical order given by $S \leq T$ if $S x \leq Tx$ for all $x \in E_+$. If $0 \leq T$, then $T$ is called positive. The dual space $E'$ equipped with the canonical order is again a Banach lattice. Instead of the operations sup and inf on $E$ we often write $\vee$ and $\wedge$, respectively. For $x \in E_+$ we denote by $[-x , x] := \{y \in E : \left| y \right| \leq x \}$ the order interval generated by $x$. A linear subspace of $E$ is an ideal if $[ -\left| x\right|, \left| x \right| ] \subseteq I$ for all $x \in I$. An ideal $I$ in $E$ is called  a band if for every subset $M \subseteq I$ such that $\sup M$ exists in $E$ one has $\sup M \in I$. A band $I$ in $E$ is called a projection band if there is a linear projection $P: E \to I$ such that $0\leq Px\leq x$ for all $x\in E_+$. Such a projection is called a band projection. A Banach lattice $E$ is called a KB-space whenever every increasing norm bounded sequence of $E_+$ is norm convergent. In particular, it follows that every KB-space has order continuous norm. All reflexive Banach lattice and AL-space are examples of KB-spaces. The following theorem is a combination of results by many mathematicians, for proofs see \cite{AB1, M, S}. 

\begin{theorem}\label{thm1}
For a Banach lattice $E$ the following statements are equivalent: 

\begin{itemize}
\item $E$ is a KB-space.
\item $E$ is a band of $E''$.
\item $E$ is weakly sequentially complete.
\item $c_0$ is not embeddable in $E$.
\item $c_0$ is not lattice embeddable in $E$.
\end{itemize}
\end{theorem} 

\section{LR-sequences}
The main tool in this section is the operator sequence on Banach space X. A family $\Theta = (T_n)_{n\in\mathbb N} \subseteq {\mathcal L}(X)$ is called an operator sequence. The sequence $\Theta$ is strongly convergent if the norm-limit $\displaystyle  \left\| \cdot\right\| - \lim_{n\to \infty} T_{n} x$ exists for each $x \in X$. A vector $x$ is called a fixed vector for the sequence $\Theta$ if $T_{n}x = x$ for each $n\in \mathbb N$. We denote by $\mathrm{Fix} (\Theta)$ the set of all fixed vectors of $\Theta$. 

The following important concept was introduced by H.P.Lotz \cite{Lo} and F. R\"abiger \cite{R2}, we use the modified terminology from \cite{EE2} called LR-sequences.

\begin{definition} \label{defn1}
A sequence $\Theta =(T_n)_{n\in\mathbb N}$ is called LR-sequence if 
\begin{itemize}
\item[1.]  $\Theta$ is uniformly bounded;
\item [2.]  $\displaystyle  \lim_{n\rightarrow\infty} \left\| T_{n} \circ (T_{m} - I) x \right\|= 0$  for every $m$ and for every $x \in X$;
\item [3.]  $\displaystyle  \lim_{n\rightarrow\infty} \left\| (T_{m} - I) \circ T_{n} x \right\|= 0$  for every $m$ and for every $x \in X$.
\end{itemize} 
\end{definition} 

 Many examples of LR-sequences appear in the investigation of operator semigroups. The Cesaro averages of a power bounded operator form an LR-sequence and moreover encompass Cesaro averages of higher orders. Also for the conditional expectation $C_n$, $(I-C_n)_{n\in\mathbb N}$ forms an LR-sequence.  For more details, we refer to \cite{ BE} -\cite{ E1}.

The following theorem is the main analytic tool in the investigation of LR-sequences. For a proof, we refer to  \cite{EE2}.

\begin{theorem}\label{thmlr}
Let $\Theta$ be an LR-sequence on a Banach space $X$. Then the following conditions are equivalent:
\begin{itemize}
\item[i.] The sequence $\Theta$ is strongly convergent. 
\item[ii.] $X = \mathrm{Fix}(\Theta) \oplus \overline{\cup_{n\in\mathbb N} (I - T_{n}) X}$ .
\item[iii.] The sequence $(T_{n} x)_{n\in\mathbb N}$ has a weak cluster point for every $x \in X$. 
\item[iv.] The fixed space $\mathrm{Fix} (\Theta)$ separates the fixed space $\mathrm{Fix} (\Theta')$ of the adjoint operator sequence $\Theta' = (T'_{n})_{n\in\mathbb N}$ in $X'$.
\end{itemize}
If one of the above vonditions holds then the strong limit of $\Theta$ is a projection onto $\mathrm{Fix}(\Theta)$. 
\end{theorem}

\section{Constrictors}
The constrictiveness of an operator is introduced in order to characterize asymptotically periodic Markov operator on $L^1$-spaces. Many authors have extended this notion to more general situations. All these notions have in common the general principal reflected by the notion of attractor introduced in \cite{LLY}.

\begin{definition}
Let $\Theta = (T_n)_{n\in\mathbb N}$ be an operator sequence on a Banach space $X$ and $C \subseteq X$. Then $C$ is called a constrictor of $\Theta$ if 
$$
\displaystyle  \lim_{n\to\infty} \mathrm{dist} (T_n  x , C) = 0
$$
for all $x \in B_X := \{z \in X : \left\|z\right\| \leq 1\}$. 
\end{definition}

Our aim is to find conditions on the constrictor $C$ implying nice asymptotic properties of $\Theta$. The first property is that every LR-sequence possessing a weakly compact constrictor is strongly convergent \cite{EE2}. Moreover, in $L^1$-spaces if $C= W + \eta B_E$ where $W$ is weakly compact and $0\leq \eta < 1$, then an LR-sequence is strongly convergent. Emel'yanov proved also that every LR-sequence containing a weakly compact operator is strongly convergent \cite{E1}.

\section{Ergodicity of LR-sequences on Banach lattices}
For our notation and terminology, we refer to \cite{AB1,M,S}. Consider, the order ideal $E_e := \cup \{[-n e, ne] : n \geq 0\}$  for any $e \in E_+$ . If $E_e$ is norm-dense in Banach lattice $E$ then $e\in E_+$ is called a quasi interior point of $E_+$. Moreover, let $\Theta$ be a positive operator sequence on $E$, then $x \in E$ is called a positive fixed vector of maximal support if $x \in \mathrm{Fix}(\Theta) \cap E_+$ and every $y \in \mathrm{Fix} (\Theta) \cap E_+$ is contained in the band generated by $x$. For every quasi-constrictive Markov operator there exists an invariant density with maximal support, see \cite{LS}. For positive contraction operator $T$ on KB-spaces with quasi-interior point, being $C:= [-z,z] + \eta B_E$ constrictor of $T$ where $z \in E_+$ and $0 \leq \eta < 1$ is proven by R\"abiger, \cite{R}. Then either $T$ is mean ergodic or there is a positive fixed vector $y \neq 0$ of $T$ of maximal support and for such positive fixed vector of maximal support $((I-P_y) A_n^T)_n $ converges strongly to zero where $P_y$ is the band projection from $E$ onto the band generated by $y$ and $A_n^T$ is the Cesaro averages of $T$.

The main theorem in \cite{EW2} is firstly generalized on  $L^1$-spaces   for Markov LR-sequences in \cite{EE1}. The principal tool in the proof of the main results of \cite{EE1} was using the additivity of the norm on the positive part of the $L^1$-space. Since this is no longer the case for a general KB-space, we use different ideas in this paper, inspried by \cite{R}. 

\begin{theorem}\label{thmc1}
Let $E$ be a KB-space with quasi-interior point $e$ and $\Theta = (T_n)_{n\in\mathbb N}$ be a positive LR-sequence in $E$, $W$ be a weakly compact subset of $E$, and $\eta \in \mathbb R$, $0\leq \eta < 1$ such that
$$
\displaystyle  \lim_{n\to\infty} \mathrm{dist}(T_n x, W + \eta B_E) =0
$$
for any $x\in B_E := \{z\in E : \left\| z\right\| \leq 1 \}$. Then $\Theta$ converges strongly. 
\end{theorem}

{\bf Proof:} The proof is motivated by the proof of Theorem 5.3 in R\"abiger's paper \cite{R}. In the first case, $(T'_{n} \phi)_{n\in\mathbb N}$ is a weak*-nullsequence for each $\phi \in E'$. Then $(T_{n} x)_{n \in \mathbb N}$ has a zero as a weak cluster point for each $x \in E$ and hence by Theorem \ref{thmlr} our $LR$-sequence converges strongly to zero. 

In the second case, there is $\phi \in E'_+$ such that $(T'_{n} \phi)_{n\in\mathbb N}$ is not $\sigma (E', E) $-convergent to $0$. Let $0 \neq \psi \in E'_+$ be a $\sigma (E', E)$-cluster point of $(T'_{n} \phi)$. We may assume $\left\| \psi \right\| = 1$. Then for all $\epsilon > 0$, there exists $n_{\epsilon}$ such that $| \left\langle \psi, x \right\rangle - \left\langle T'_{n} \phi, x\right\rangle | < \epsilon$ and $| \left\langle T'_{m} \psi, x \right\rangle - \left\langle T'_{m} T'_{n} \phi, x\right\rangle | < \epsilon$ for every $m\in\mathbb N$ and for every $n\geq n_{\epsilon}$. Therefore we get $T'_{m} \psi = \psi$. 

Now for fix $\epsilon > 0$ satisfying $\epsilon \leq 1- \eta$ choose $x \in B_E \cap E_+$ such that $\left\langle \psi, x\right\rangle > 1 - \epsilon$. Let $x'' \in E''_+$ be a weak-cluster point of $(T_{n} x)$. Then there exists $n_{\epsilon}$ such that $\left\langle \psi, x'' \right\rangle - \left\langle \psi, T_{n} x\right\rangle < \epsilon$ and $\left\langle T'_{m} \psi, x'' \right\rangle - \left\langle T'_{m} \psi, T_{n} x\right\rangle < \epsilon$ by combining these two estimates with the property of LR-sequence, we obtain $T''_{m} x'' = x''$ for every $m\in\mathbb N$. Since $\displaystyle \lim_{n\to\infty} \mathrm{dist} (T_{n} x, W + \eta B_E) =0$, $W$ is weakly compact and $E$ is a band in $E''$ means $E$ is weak$^*$-closed in $E''$ then we obtain $x'' \in W + \eta B_E$. Moreover $x''$ is a weak*-cluster point of $(T_{n} x)$, then for every $\epsilon' > 0$, there exists $n_{\epsilon}$ such that $| \left\langle \psi, x'' \right\rangle - \left\langle \psi, T_{n_{\epsilon}} x\right\rangle | < \epsilon'$. Therefore we have $| \left\langle \psi, x'' \right\rangle - \left\langle T'_{n_{\epsilon}} \psi, x\right\rangle | < \epsilon'$ and since $T'_{n} \psi = \psi$, we obtain $| \left\langle \psi, x'' \right\rangle - \left\langle \psi, x\right\rangle | < \epsilon'$. By arbitrariness of $\epsilon'$, $\left\langle \psi, x'' \right\rangle = \left\langle \psi, x\right\rangle $.

Being $E$ a KB-space, by Theorem \ref{thm1} $E$ is a band in $E''$. Denote by $P$ the band projection from $E''$ onto $E$  i.e. $P : E'' \to E$. Hence 
\begin{eqnarray} \label{1}
\left\langle \psi, Px'' \right\rangle &=& \left\langle \psi, x'' \right\rangle - \left\langle \psi, (I_{E''} - P)x'' \right\rangle \\
&=& \left\langle \psi, x \right\rangle - \left\langle \psi, (I_{E''} - P) x'' \right\rangle \nonumber \\
&>& 1- \epsilon - \eta > 0 \nonumber
\end{eqnarray}

It follows from (\ref{1}) that $P x'' \neq 0$. Since $x''$ is a weak*-cluster point of $(T_{n} x)_{n\in\mathbb N}$, $P x'' > 0$ and moreover, since $E$ has order continuous norm $\displaystyle  z := \lim T_{n} P x'' \in E_+$ exists. Clearly $T_{m} z = z$ and from $\left\langle \psi, z \right\rangle = \left\langle \psi, Px'' \right\rangle > 0$, it follows that $z \neq 0$. Hence $\mathrm{Fix} (\Theta) \cap E_+ \neq \{0\}$. 

The existence of a quasi-interior point $e$ implies the existence of a strictly positive linear functional $\psi$ [\cite{AB1}, Theorem 12.14]. For $x\in E$, let $P_x$ be the band projection from $E$ onto the band generated by $x$. Consider the element $\displaystyle \alpha := \sup_{x \in Fix(\Theta) \cap E_+} \left\langle \psi, P_x e \right\rangle > 0$. Choose a sequence $x_{n} \in \mathrm{Fix}(\Theta) \cap E_+$, $n\in \mathbb N$, $\left\| x_{n} \right\| \leq 1$ and $\displaystyle  \alpha = \lim \left\langle \psi, P_{x_{n}} e \right\rangle$. Define $u = \sum_{n} 2^{-n} x_{\lambda_n}$, then $u$ is also an element of $\mathrm{Fix}(\Theta) \cap E_+$ and in addition $P_u \geq P_{x_{n}}$ for all $n\in\mathbb N$. Hence $\left\langle \psi, P_u e \right\rangle = \alpha$. 

Further taking $x \in \mathrm{Fix}(\Theta) \cap E_+$, clearly $P_{u+x} \geq P_x$ and $P_{u+x} \geq P_u$. From the limit property of $\left\langle \psi, P_{x_{n}} \right\rangle$, $\alpha \leq \left\langle \psi, P_{u+x} e\right\rangle \leq \alpha$, so $\alpha = \left\langle \psi, P_{u+x} e\right\rangle$ and we know above $\psi$ is strictly positive, it implies that $P_u e = P_{u+x} e$. Owing to quasi-interior point $e$, $P_u = P_{u+x}$ and by $P_{u+x} \geq P_x$ then we obtain $P_u \geq P_x$. Hence $u$ has a maximal support. 

In the next step, we will prove that for the band projection denoted by $P_u$, of positive fixed vector of maximal support, $((I- P_u) T_n)$ converges to zero strongly as $n\to\infty$. Let $P_u$ be a band projection onto $B_u$ where $B_u = \overline{\cup_n [-n u, nu]}$. Denote the new operator $Q = I_E- P_u$ and the sequence ${\mathcal S} = (S_{n}) = (Q T_{n})$. Since $u$ is a fixed vector so $T_{n} B_u \subset B_u$ and hence we get $T_{n} P_u = P_u T_{n} P_u$ and in addition $Q T_{n} Q = Q T_{n}$ for each $n$. 

Our aim is to show that $(QT_{n}) = (S_{n})$ converges strongly to zero. If not, then there exists by above in the second case of proof, $0 \neq \psi \in \mathrm{Fix}{{\mathcal S}'} \cap E'_+$. 
$$
\psi = S'_{n} \psi = T'_{n} Q' \psi = Q' T'_{n} Q' \psi = Q' S'_{n} \psi = Q' \psi
$$
and hence $\psi = S'_{n} \psi = T'_{n} Q'\psi = T'_{n} \psi$, namely, $\psi \in \mathrm{Fix}(\Theta')$. By above, there is $0 \neq x \in \mathrm{Fix} (\Theta) \cap E_+$ such that $\left\langle x, \psi \right\rangle > 0 $. Then
$$
0 < \left\langle x, \psi \right\rangle = \left\langle x, Q' \psi \right\rangle = \left\langle Qx, \psi \right\rangle 
$$
implies $Q x \neq 0$, i.e. $x\notin B_u$. It is a contradiction to our assumption on $u$. Indeed, $S_{n} \to 0$ strongly as $n\to\infty$.

In the next step, we will prove that $(T_n)$ converges strongly by using Theorem \ref{thmlr}.

We know that our operator sequence is an LR-sequence and positive. For fixed $\epsilon > 0$ and $x\in E$, $\displaystyle  \lim_{n\to \infty} S_{n} x = \lim_{n\to\infty} (I-P_u) T_{n} x = 0$, there exists $n_{\epsilon}$ such that $\mathrm{dist} (T_{n_{\epsilon}} x, B_u) \leq \frac{\epsilon}{3 M}$ where $\displaystyle M = \sup_{n} \left\| T_{n} \right\|$. It implies that there exists $c_{\epsilon} \in {\mathbb R}_+$ and $y \in [- c_{\epsilon} u, c_{\epsilon} u]$ satisfying $\left\| T_{n_{\epsilon}} x - y \right\| \leq \frac{\epsilon}{2M}$. 

For any $m\in\mathbb N$, $\left\| T_{m} T_{n_{\epsilon}} x - T_{m} y \right\| \leq \left\| T_{m} \right\| \left\| T_{n_{\epsilon}} x- y \right\| \leq \frac{\epsilon}{2} $. Moreover, since $[-u,u]$ is $\Theta$-invariant then we get $T_{n} y \in [-c_{\epsilon}u, c_{\epsilon} u]$ for each $n$. Therefore $$\mathrm{dist} (T_{n_{\epsilon}} x, [-c_{\epsilon} u, c_{\epsilon} u]) \leq \epsilon'$$, i.e. for any $\epsilon' > 0 $, there exists an interval $[-c_{\epsilon}, c_{\epsilon}]$ such that $(T_{n} x)_{n=0}^{\infty} \subseteq [-c_{\epsilon}, c_{\epsilon}] + \epsilon' B_E$. It shows that $(T_{n} x)$ has a weak cluster point because $E$ is a KB-space and almost order bounded subset of $E$ is weakly precompact. Then by Theorem \ref{thmlr}, $(T_{n}x)_{n\in\mathbb N}$ is norm convergent for any $x\in E$, i.e. $\Theta$ converges strongly.
\begin{flushright}
$\square$
\end{flushright}

The previous theorem can be formulated for the net case as follows:
\begin{theorem}
Let $E$ be a KB-space with quasi-interior point $e$ and $\Theta = (T_{\lambda})_{\lambda\in\Lambda}$ be a positive LR-net in $E$ which has a cofinal subsequence, $W$ be a weakly compact subset of $E$, and $\eta \in \mathbb R$, $0\leq \eta < 1$ such that
$$
\displaystyle  \lim_{\lambda\to\infty} \mathrm{dist}(T_{\lambda} x, W + \eta B_E) =0
$$
for any $x\in B_E := \{z\in E : \left\| z\right\| \leq 1 \}$. Then $\Theta$ converges strongly. 
\end{theorem}

The theorem is also true if we replace a weakly compact subset $W$ of $E$ by an order interval $[-g, g]$ for any $g \in E_+$ because in KB-spaces, every order intervals are weakly compact. Besides in this case we have more results that also dimension of fixed space is finite. 

\begin{theorem} \label{thmc2}
Let $E$ be a KB-space with a quasi-interior point $e$, $\Theta=(T_n)_{n\in\mathbb N}$ be a positive LR-sequence. Then the following are equivalent 
\begin{itemize}
\item[i] there exists a function $g \in E_+$ and $\eta \in \mathbb R$, $0\leq\eta < 1$ such that
$$
\displaystyle  \lim_{n \to\infty} \mathrm{dist}(T_{n} x, [-g, g] + \eta B_E) =0, \,\,\, \forall x \in B_E
$$  
\item[ii] the sequence $\Theta$ is strongly convergent and $\mathrm{dim Fix}(\Theta) < \infty$.
\end{itemize}
\end{theorem}

{\bf Proof:} 
The proof of this theorem is motivated by the proof of Theorem 3 in \cite{EE1}.

$(i) \Rightarrow (ii) :$ By the previous theorem, Theorem \ref{thmc1},  $\Theta$ converges strongly onto  $\mathrm{Fix}(\Theta)$.
Therefore for each $x \in B_E$ $Px \in [-g,g] + \eta B_E$ by the statement of the theorem. If we consider the iteration of $Px$ then we obtain
$$
\displaystyle
Px = P^2 x  \in [-Pg, Pg] + \eta P(B_E).
$$
Since $P$ is a projection, i.e. $\left\| P\right\| \leq 1$,
$$
\displaystyle
Px = P^2 x \in [-Pg, Pg] + [-\eta g, \eta g] + \eta^2 B_E.
$$
 If we repeat of the iteration, we have for arbitrary $n\in\mathbb N$,
$$
\displaystyle
Px = P^n x \in [-\sum_{i=0}^{n-1} \eta^{i} P^{n-i} g, \sum_{i=0}^{n-1} \eta^{i} P^{n-i} g] + [- \eta^{n-1} g, \eta^{n-1} g] + \eta^n B_E.
$$
Hence
$$
\displaystyle
Px \in [-\sum_{i=0}^{n-1} \eta^{i} P g, \sum_{i=0}^{n-1} \eta^{i} P g] + [- \eta^{n-1} g, \eta^{n-1} g] + \eta^n B_E.
$$
 If we continue to iterate $Px$ we get the condition that $Px \in [-c,c]$ where $\displaystyle c=\frac{1}{1-\eta} Pg$. Hence $P(B_E) \subseteq [-c,c]$ that is to say $\mathrm{Fix}(\Theta)$ is finite dimensional space, \cite{S}.

$(ii) \Rightarrow (i):$ If $\mathrm{dim (Fix} (\Theta)) < \infty$, then there exists a family of pairwise disjoint densities $u_1, u_2, \cdots, u_k$ such that $\mathrm{Fix}(\Theta) = span\{u_1,u_2,\cdots,u_k\}$. Denote the element $g := u_1+\cdots+u_k$ and taking an element from $B_E \cap E_+$, then $\displaystyle  P x := \lim T_n x$ is a linear combination of  $u_1, \cdots , u_k$ say $Px = \sum_{i=1}^k \alpha_i u_i \leq \sum_{i=1}^k u_i$. Thus 
$$
\displaystyle  \limsup_{n \to\infty} \left\| (T_n x - g)_+ \right\| = \left\| (Px - g)_+\right\| = 0
$$
for every $x \in B_E \cap E_+$.  
\begin{flushright}
$\square$
\end{flushright}

In the above two theorems, KB-space conditions cannot be omitted. Even for Banach lattices with order continuous norm, this result can fail, for the counterexample, see \cite{EW1}.

Additionally we also can be formulated the above theorem for the net case as follows:
\begin{theorem}
Let $E$ be a KB-space with quasi-interior point $e$ and $\Theta = (T_{\lambda})_{\lambda\in\Lambda}$ be a positive LR-net in $E$ which has a cofinal subsequence. Then the following are equivalent 
\begin{itemize}
\item[i] there exists a function $g \in E_+$ and $\eta \in \mathbb R$, $0\leq\eta < 1$ such that
$$
\displaystyle  \lim_{\lambda \to\infty} \mathrm{dist}(T_{\lambda} x, [-g, g] + \eta B_E) =0, \,\,\, \forall x \in B_E
$$  
\item[ii] the sequence $\Theta$ is strongly convergent and $\mathrm{dim Fix}(\Theta) < \infty$.
\end{itemize}
\end{theorem}

\section{Asymptotic stability of LR-sequences}
The asymptotic stability of positive operators and lower-bound technique is developed  in applications of Markov operators. In this section we prove the following theorems as a corollary of Theorem \ref{thmc2}. Theorem \ref{thmc3} is the generalization of Theorem 4 in \cite{EE1}. Emelyanov and Erkursun proved asymptotic stability and existence of lower bound function are equivalent for Markov LR-nets on $L^1$-spaces. In this section we will have a KB-space as well. 

In the first, we give the following two definitions which are motivated for operator nets by the definitions used in \cite{LM}.

\begin{definition}
Let $\Theta$ be an operator nets on KB-spaces.  $\Theta$ is called asymptotically stable whenever there exists an element $u \in E_+ \cap U_E$ where $U_E := \{x \in E : \left\| x\right\| = 1 \}$ such that
 $$
 \lim_{\lambda \to \infty} \left\| T_{\lambda} x  - u \right\| = 0
 $$
 for every element from $E_+ \cap U_E$.
\end{definition}

\begin{definition}
An element $h \in E_+$ is called lower bound element for $\Theta$ if 
 $$
 \lim_{\lambda\to\infty} \left\| (h - T_{\lambda} x)_+ \right\| =0
 $$
for every element $x \in E_+ \cap U_E$
\end{definition}

For main results, positivity is not only sufficient in addition we need Markov operators. Before proving the theorem, we need to define Markov operator nets on a Banach lattice $E$.

\begin{definition}
Let $E$ be a Banach lattice.  A positive, linear, uniform bounded operator net  $\Theta =(T_{\lambda})_{\lambda\in\Lambda}$
is called a Markov operator net if there exists a strictly positive element  $0 < e' \in  E'_+$ such that $T_{\lambda}' e' = e'$ for each $\lambda\in\Lambda$.
\end{definition}

If we consider in the sequence case $\Theta = (T_n)_{n\in\mathbb N}$, each element of $\Theta$ is Markov operator, still we need a common fixed point $e'$. For instance even on $L_1$-space, the element $T_m$ of $\Theta$ might not be a Markov operator on the new norm space $(L_1, e_n')$ for each $n\neq m \in \mathbb N$.

As a remark if we consider the net as $\Theta = (T^n)_{n\in \mathbb N}$ as the iteration of a single operator $T$, Markov operator sequence $\Theta$ means $T$  is power-bounded. It is the general version of the definition in \cite{H} where $T$ is contraction.
It is well known that if $T$ is a positive linear operator defined on a Banach lattice $E$, then $T$ is continuous. It is also well known that if the Banach lattice $E$ has order continuous norm, then the positive operator $T$ is also order continuous. We note that the Markov operators, according to this definition, are again contained in the class of all positive power-bounded and that the adjoint $T'$ is also a positive and power-bounded. For more details, we refer to \cite{H}.

In the following theorem, we will establish the asymptotic properties of Markov LR-sequences. It is firstly given on $L^1$-spaces in \cite{EE1} which are the examples of Markov LR-sequences which need not to be $\mathcal T$-ergodic sequences.  Now we will prove this results on KB-spaces. Before them, we need technical tools for proof. The technical lemma connects norm convergence of order bounded sequences in KB-spaces with convergence in $(E, e')$ for suitable linear forms $e' \in E'$. Recall that $e' \in E'$ is strictly positive if $\left\langle x, e'\right\rangle >0$ for all $x\in E_+ \setminus \{0\}$. We refer to \cite{R} for proof of the lemma.

\begin{lemma} \label{lemma}
Let $(x_{n})_{n\in\mathbb N}$ be an order bounded sequence in a KB-space and let $x' \in E'$ be strictly positive. Then $\displaystyle  \lim_{n\to\infty}\left\| x_n \right\|=0 $ if and only if $\displaystyle \lim_{n\to\infty} \left\langle \left| x_n \right| , x' \right\rangle =0$.
\end{lemma}

\begin{theorem} \label{thmc3}
Let $\Theta$ be a Markov LR-sequence on KB-spaces with fixed common element $e'$. Then the following are equivalent:
\begin{itemize}
\item[i] $\Theta$ is asymptotically stable
\item[ii] $\Theta$ has nontrivial lower-bound element in the space $(E, e')$
\item[ii] $\Theta$ has nontrivial lower-bound element in the space $(E, e')$ for each $e' \in E'_+$.
\end{itemize}
\end{theorem}

{\bf Proof:} 
$(i) \Rightarrow (iii) :$ Let $g \in E_+ \cap U_E$  satisfy
 $$
 \lim_{t \to \infty} \left\| A^T_t x  - g \right\| = 0
 $$
 for every $x \in E_+ \cap U_E$, then $g$ is automatically a nontrivial lower-bound function for $\Theta$ on the space $(E,e')$ for each strictly positive $e' \in E'_+$.

$(iii) \Rightarrow (ii) :$ Obvious

$(ii) \Rightarrow (i) :$ Let $h$ be a lower-bound element of $\Theta$ in $(E , e')$, i.e., 
$$
\displaystyle \lim_{n\to\infty}  \left\langle (T_n x - h)_+ , e' \right\rangle =0 \,\,\,\,\, \forall x\in E_+ \cap U_E. 
$$
Since the norm on $(E,e')$ is an $L_1$-norm, then we can consider
$$
\displaystyle \limsup_{n\to\infty}  \left\langle (T_n x - h)_+ , e' \right\rangle \leq \eta
$$
 where $\eta :=1- \left\| h \right\|_{(E,e')} = 1 -  \left\langle h , e' \right\rangle$ for each $x \in E_+ \cap U_E$.

 By Theorem \ref{thmc2}, ${\tilde T}_n$ where $\tilde{T} j_{e'} x = j_{e'} Tx $ for lattice homomorphism $j_{e'} : E \to (E, e')$  converges strongly to the finite dimensional fixed space of $\tilde {\Theta}$. Therefore by Eberlein's Theorem 
$$
(E, e') = \mathrm{(Fix}(\tilde {\Theta})) \oplus \mathrm{Ker}(\tilde {\Theta}).
$$
In addition we know that $\mathrm{Fix}(\tilde {\Theta})$ is a sublattice of $(E,e')$ and by Judin's Theorem, it possesses a linear basis $\tilde{(u_i)}_{i=1}^n$ where $n = \mathrm{dim (Fix} (\tilde {\Theta}))$ which consists of pairwise disjoint element with $\left\|\tilde{ u_i }\right\| =1$, $i=1, \cdots, n$, see \cite{S}. Since $T u_i = u_i$ for each $i=1,\cdots, n$, 
$$
 \left\langle (h - u_i)_+,   e' \right\rangle = \left\langle (h- T u_i)_+,  e' \right\rangle = \lim_{t\to\infty}  \left\langle (h - T_n  u_i)_+ ,  e' \right\rangle =0
 $$
 implies 
\begin{eqnarray} \label{eqn1}
u_i \geq h \geq 0 \,\,\,\,\,\,\,\, i =1, \cdots, n.
\end{eqnarray}
Since $\tilde{(u_i)}_{i=1}^n$ is pairwise disjoint with $\left\| \tilde{u_i} \right\|_{(E,e')} =1$ the condition \ref{eqn1} ensure that $\mathrm{dim (Fix}(\tilde {\Theta})) =1$. Therefore $\mathrm{Fix} (\tilde {\Theta}) = \mathbb{R}\tilde{ u_1}$ and for every element $x \in E_+ \cap B_E$, $\lim_{n\to\infty} T_n x = u_1$.
\begin{flushright}
$\square$
\end{flushright}

The Lasota's criterion of asymptotic stability says that a one-parameter Markov semigroup if and only if there is a nontrivial lower-bound function. In \cite{EE1} Lasota's lower-bound criteria is generalized on $L^1$-spaces to abelian Markov semigroups. In this proposition we generalize it on KB-spaces. An abelian Markov semigroup is an operator net with respect to the natural partial order $\succ$ mentioned Section 2.

\begin{proposition} \label{prop1}
Every asymptotically stable abelian Markov semigroup ${\mathcal T} = (T_t)_{t}$which has a common fixed point $e'$ of ${\mathcal T}'$ on KB-space is an LR-net.
\end{proposition}

{\bf Proof:}
Since any Markov operator on KB-space is a positive contraction, a Markov sequence is uniformly bounded. We need to check (LR2) and (LR3) conditions of Definition \ref{defn1}. Moreover because of abelian property, it suffices to prove only (LR2) or (LR3).  Without loss of generality, taken for an arbitrary element $x$ from $E_+ \cap U_E$. Since $\mathcal T$ is asymptotically stable then by Lemma \ref{lemma}  $\displaystyle \lim_{n\to\infty} \left\langle  |T_{t_n} x - u|, e' \right\rangle =0$. If we consider fixed $t_m\in\mathbb R$, $ \left\langle ((I- T_{t_m}) x)_+ , e' \right\rangle = \left\langle ((I-T_{t_m})x)_- , e' \right\rangle $ by additivity property of $(E,e')$.

Now for (LR2), call $(I-T_{t_m})x =y$
\begin{eqnarray*}
\left\langle |T_{t_n} (I - T_{t_m}) x|, e' \right\rangle &=& \left\langle |T_{t_n} y|, e' \right\rangle = \left\langle |T_{t_n} y_+ - T_{t_n} y_-|, e' \right\rangle \\
&=& \left\langle | T_{t_n} y_+ - \left\| y_+ \right\|_{(E,e')} u  + \left\| y_- \right\|_{(E,e')} u  -T_{t_n} g_- | , e' \right\rangle \\
&\leq& \left\langle | T_{t_n} y_+ - \left\| y_+ \right\|_{(E,e')} u | , e' \right\rangle \\
  & & + \left\langle | T_{t_n} y_- - \left\| y_-\right\|_{(E,e')} u | , e' \right\rangle 
\end{eqnarray*}
which converge to zero as $n \to \infty$. Therefore $\lim_{n\to\infty} T_{t_n} (I- T_{t_m}) x =0$ for each $x \in U_E \cap E_+$.
\begin{flushright}
$\square$
\end{flushright}

\begin{proposition} \label{prop2}
Let ${\mathcal T} = (T_t)_{t}$ be an abelian Markov semigroup which has a common fixed point $e'$ of ${\mathcal T}'$ on KB-spaces possessing a nontrivial lower-bound function on the space $(E,e')$, then $\mathcal T$ is an LR-net.
\end{proposition}

{\bf Proof:}
The important part of the proof is that since $T_{t'}$ is Markov then $\left\langle ((I-T_{t'}) f)_+, e' \right\rangle = \left\langle ((I-T_{t'}) f)_-, e' \right\rangle$ for each $t'$ and for each $f\in E$ on the space $(E,e')$. It implies that $\left\| ((I-T_{t'}) f)_+ \right\|_{(E,e')} = \left\| ((I-T_{t'}) f)_- \right\|_{(E,e')}$. We repeat the argument in \cite{EE1} in short for convenient of the reader.

Let $0 \neq h \in E_+$ be a nontrivial lower-bound element for $\mathcal T$ on the space $(E,e')$, then $\left\| h\right\|_{(E,e')} \leq \left\| (h -T_t f)_+ \right\|_{(E,e')} + \left\| h \wedge T_{t} f \right\|_{(E,e')} \leq \epsilon + 1 $ for every $t$ so obviously $\left\| h\right\|_{(E,e')} \leq 1$. Since semigroup is Markovian then it is uniformly bounded. Moreover because of abelian property, it suffices to prove only (LR2) or (LR3). Thus we prove the following formula.
\begin{eqnarray} \label{eqns1}
\lim_{t\to\infty} \left\| T_t (I - T_{t'}) f \right\|_{(E,e')} = 0 \,\,\,\,\,\,\, (\forall t', f \in E)
\end{eqnarray}

Take any element $f \in B_E$, then we know that $T_{t'}$ is Markov and then
\begin{eqnarray} \label{eqns2}
\left\langle ((I-T_{t'}) f)_+, e' \right\rangle = \left\langle ((I-T_{t'}) f)_-, e' \right\rangle
\end{eqnarray}

Therefore by \ref{eqns1}, we have to prove that $\lim{t\to\infty} \left\| T_t f\right\|_{(E,e')} =0$ for every $f \in B_E$ such that \ref{eqns2} holds. Define the set $E_0 := \{ f\in E : \left\| f_+ \right\|_{(E,e')}=\left\| f_-\right\|_{(E,e')} \}$.

Take any element $f\in E_0$ such that $f =2^{-1} \left\| f\right\|_{(E,e')} (f_1 - f_2)$ where $f_1 = 2 \left\| f\right\|_{(E,e')}^{-1} f_+$ and $f_2 = 2 \left\| f\right\|_{(E,e')}^{-1} f_-$. Hence $f_1$ and $f_2$ are elements of $E_+ \cap U_E$.

Since $h$ is the lower-bound element for the Markov semigroup $\mathcal T$, there exists $t_1$ such that $\left\| (h - T_{t_1} f_1)_+  \right\|_{(E,e')} \leq \frac{1}{4} \left\| h \right\|_{(E,e')}$ and $\left\| (h - T_{t_1} f_2)_+ \right\|_{(E,e')}\leq \frac{1}{4} \left\| h \right\|_{(E,e')}$ hold for every $t \geq t_1$. From Riesz space properties, we obtain $\left\| T_{t_1} f_1- T_{t_1} f_2\right\|_{(E,e')} \leq 2- \frac{1}{2} \left\| h\right\|_{(E,e')}$ and $\left\| T_t f \right\|_{(E,e')} \leq (1 - \frac{1}{4} \left\| h\right\|_{(E,e')}) \left\| f\right\|_{(E,e')} $ for every $t \geq t_1$.

Replacing $f$ with $T_{t_1} f$ which is also an element of $E_0$ and repeating the argument above gives an element $t_2$ such that 
$$
\left\| T_t T_{t_1} f\right\|_{(E,e')} \leq (1 - \frac{1}{4} \left\| h\right\|_{(E,e')}) \left\| T_{t_1} f \right\|_{(E,e')} \,\,\,\,\, \forall t \geq t_2
$$

By induction, we can generate a sequence $(t_n)$ such that
\begin{eqnarray*}
\left\| T_{t} f \right\|_{(E,e')} \leq \left\| T_t T_{t_{n-1}} f\right\|_{(E,e')} &\leq& (1- \frac{1}{4} \left\| h\right\|_{(E,e')}) \left\| T_{t_{n-1}} f \right\|_{(E,e')} \\
&& \vdots \\
 (\forall t \geq t_1 + \cdots + t_n) \, &\leq& (1-\frac{1}{4} \left\| h\right\|_{(E,e')})^n  \left\| f\right\|_{(E,e')} 
\end{eqnarray*}

Since $0< \left\| h\right\| <1$, then $\lim_{t\to\infty} \left\| T_t f\right\|_{(E,e')} =0$ and hence the proof is completed.

\begin{flushright}
$\square$
\end{flushright}

\begin{theorem} \label{thmc4}
Let ${\mathcal T} = (T_t)_t$ be an abelian Markov semigroup which has a common fixed point $e'$ of ${\mathcal T}'$ on KB-spaces $E$.  Then the following are equivalent:
\begin{itemize}
\item[i] $\mathcal T$ is asymptotically stable
\item[ii] There exists a nontrivial lower-bound element for $\mathcal T$ in the space $(E,e')$
\end{itemize}
\end{theorem}

{\bf Proof:}
Since the asymptotic stable Markov semigroup $\mathcal T$ is the LR-sequence by \ref{prop1}, the existence of nontrivial lower-bound element for $\mathcal T$ follows from Theorem \ref{thmc3}. In addition, the existence of nontrivial lower bound element for $\mathcal T$ gives us that $\mathcal T$ is an LR-sequence by Proposition \ref{prop2} and the asmyptotic stability of $\mathcal T$ follows from Theorem \ref{thmc3} .

\begin{flushright}
$\square$
\end{flushright}

{\bf Acknowledgment:} The author gratefully acknowledge Prof. Eduard Emelyanov for his helpful remarks.

\end{document}